\documentclass[12pt, twoside]{article}
 \usepackage{amsfonts, amssymb, amsmath, amsthm, latexsym, amscd}
\theoremstyle{plain}
\newtheorem{Th}{Theorem}[section]
\newtheorem{Lem}[Th]{Lemma}

\newtheorem{Cor}[Th]{Corollary}
\theoremstyle{definition}

\begin{document}
\addtolength{\textheight}{0 cm} \addtolength{\hoffset}{0 cm}
\addtolength{\textwidth}{0 cm} \addtolength{\voffset}{0 cm}

\title{Existence and nonexistence of solutions for a singular $p$-Laplacian
Dirichlet problem }


\def\dis{\displaystyle}

\author{ Mahmoud Hesaaraki \textsc \bf \footnote{The author  would like to thank the Sharif University of Technology
for supporting this research.} \\
\small Department of Mathematics,  \\
\small Sharif University of Technology  \\
\small P.O. Box 11365-9415, Tehran, Iran \\
\tt hesaraki@sina.sharif.ac.ir \\
Abbas Moameni \\
 \small Department of Mathematics,  \\
\small Sharif University of Technology  \\
\small P.O. Box 11365-9415, Tehran, Iran \\
\tt moameni@math.ubc.ca}
\date{}

\maketitle

\renewcommand{\theequation}{\arabic{section}.\arabic{equation}}

\begin{abstract}
 We study the existence of positive
radially symmetric solution for the singular $p$-Laplacian Dirichlet
problem, $-\bigtriangleup_p u =\lambda |u|^{p-2} u-\gamma
u^{-\alpha}$ where $\lambda>0,\gamma>0$ and, $0<\alpha<1$, are
parameters and $\Omega$, the domain of the equation, is a ball in
$\mathbb{R}^N$. By using some variational methods we show that, if
$\lambda$ is contained in some interval, then the problem has a
radially symmetric positive solution on the ball. Moreover, we
obtain a nonexistence result, whenever $\lambda \leq 0, \gamma<0$
and $\Omega$ is a bounded domain, with smooth boundary.
\end{abstract}


\section{Introduction}
In this paper we study the singular $p$-Laplacian Dirichlet problem
\begin{equation}
\label{1.1}
\begin{array}{rl}
-\bigtriangleup_p u & = \lambda |u|^{p-2} u- g(u) \ \text{in}\
\Omega , \\
u & =0 \ \text{on}\ \partial\Omega , \\
\end{array}
\end{equation}
where $\Omega$ is a ball with center 0 in $\mathbb{R}^N$, $N\geq
2$ and $g:(0,\infty) \to (0,\infty)$ is a function satisfying
$g(\tau) \to \infty $ as $ \tau \to 0$.

 Indeed, we obtain
 existence and nonexistence results under some assumptions on
$N,p,g,\lambda$ and $\Omega$.
 Chen  in \cite{1}, in the case $p=2$ and
$g(\tau) = {\tau^{-\alpha} \over 1+\alpha}$ for $\tau>0$ and
$\Omega=\{x\in \mathbb{R}^N : |x|<R\}$, by using the shooting method
obtained the following results:
 \begin{itemize}
 \item
 { \it There are real numbers
$R_1>R_2>0$ such that the problem (\ref{1.1})  has a radially
symmetric, positive solution if $R_1\geq R>R_2$. Besides, if $u$
is a radially symmetric, positive solution for the problem in the
case of $R=R_1$, then ${\partial u\over
\partial r}=0$ on $\partial \Omega$, where ${\partial \over
\partial r}$ is the  outward normal derivative.}
\end{itemize}

In order to show the existence of solutions, we use the
variational methods by considering the following functional:
\begin{equation}
 F(u)= {1\over p} \int_{\Omega} |\bigtriangledown u|^p dx-
{\lambda \over p} \int_{\Omega} |u|^p dx+\int_{\Omega}
\int_0^{u(x)} g(\tau) dt dx, u\in W_0^{1,p} (\Omega),
\end{equation}
 associated with the problem. Since this functional is not even G\^ateaux
 differentiable, we cannot use the
 deformation argument. Neither can we use the strong  maximum
 principle because  of the property of the nonlinear term $g$.
 Here, we will show that if $u$ is a function which is a
  minimax value of $F$, then $u$ is a radially
 symmetric, positive solution of the problem.

 For the nonexistence result we use the Pohozaev identity which is
  introduced in \cite{3} and we show that if $\lambda \leq 0$, we may
 have no positive solution in $W^{2, p_1} (\Omega)\cap W_0^{1, p_1}(\Omega)$, $p_1>N$. In this case we assume that
  $\Omega$ is a bounded domain and its boundary, $\partial \Omega$
 has the following property:\\
  There  exists a unit normal vector
 $v(x)= (v_1(x), \ldots, v_N(x))$  at every point $x\in \partial
 \Omega$ such that $\sum\limits_{i=1}^N x_i v_i (x) \geq 0$.

\section{Existence Result}
In this section we prove the existence of radially symmetric
positive  solution for the problem (\ref{1.1}) in the following
theorem.

\begin{Th} Let $\Omega$ be a ball in $\mathbb{R}^N$ with center 0,
$0<\alpha<1$, $p\geq 2 $ and  $N\geq 2$. Suppose $g$ is a
 $c^{\infty}(0,\infty)$ function with $g(\tau)>0$ and
 \begin{equation}
{d\over d\tau} (g(\tau))<0,
\end{equation}
 moreover
\begin{equation}
 m_1 \tau^{-\alpha} \leq g (\tau) \leq m_2 \tau^{-\alpha}
 \end{equation}
 for some  positive constants $m_1$ and $m_2$ with ${p m_1 \over
 1-\alpha}> m_2 \geq m_1$. If
 \begin{equation}
 \lambda- { \lambda p (1-\alpha) m_1 \over p N m_2-
 (N-p)(1-\alpha) m_1} \leq \lambda_1 <\lambda ,
 \end{equation}
 where $\lambda_1$ is the first eigenvalue
 of the operator $-\bigtriangleup_p$ with homogenous Dirichlet
 boundary condition, then problem (\ref{1.1}) has a radically symmetric positive
 solution.
 \end{Th}

In order to the proof this Theorem we need some preliminary lemmas.
The following sets  will be used in our proofs.

 \[ U=\{ u \in W_0^{1,p} (\Omega): u\ \text{is\ radially\
 symmetric} \}, \]

 and
 \[ W=\{ u\in U: \int_{\Omega} |\bigtriangledown u|^p dx <\lambda
 \int|u|^p dx\}. \]
 Note that for a function $u\in W$, we may regard it as a
 one  variable function $u(r)$, where $r= |x|$ with $x\in
 \Omega$. Also
note  that $W$ is not empty, since $\lambda_1 <\lambda$ and the
eigenfunctions of $-\bigtriangleup_p$ with homogeneous Dirichlet
 boundary condition for $\lambda_1$ are radially symmetric.

In the following Lemmas we assume that all of the conditions of
Theorem 2.1 hold. Moreover, we assume that $g$ is defined on
$\Bbb{R}$ with $g(0)=0$ and $g(t)=-g(-t)$ for $t<0$.

\begin{Lem}  Let $u\in W$, then $\int_{\Omega} g(t u)
 {u\over t^{p-1}} dx \to +\infty$ as $t\to 0^+$,
 $ \int_{\Omega} {g(tu)u \over t^{p-1}} dx \to 0$ as
 $t \to +\infty$ and the function $t\to \int_{\Omega} {g
 (tu)u\over t} dx$ is strictly decreasing for $t>0$.  Especially,
 there exists a unique $t>0$ such that
 \[ \int_{\Omega} |\bigtriangledown t u|^p dx+ \int_{\Omega} g(tu)
 tu dx= \lambda \int |tu|^p dx, \]
 which is equivalent to $F(tu) =\max\limits_{s>0} F(su)$.
 \end{Lem}

 \paragraph{Proof.} From (2.4), we have
 \[ m_1 t^{-(p-1)-\alpha}  \int_{\Omega} |u|^{1-\alpha} dx \leq
 \int_{\Omega} g (tu) {u\over t^{p-1}} dx \leq m_2 t^{-
 (p-1)-\alpha} \int_{\Omega} |u|^{1-\alpha} dx, \]
for every $t>0$. Thus  we obtain
\[ \int_{\Omega} g(tu) {u\over t^{p-1}} dx \longrightarrow +\infty
\quad \text{as}\ t \longrightarrow 0^+, \]
\[ \int_{\Omega} g(tu) {u\over t^{p-1}} dx \longrightarrow 0
\quad \text{as}\ t \longrightarrow \infty.  \]
 From  (2.3), we
 see that the function $t \to \int_{\Omega}
  g(tu) {u\over t^{p-1}} dx$ is strictly decreasing for $t>0$.
  $\Box$

  We define a subset $V$ of $W$ by
  \[ V= \{u\in W: \int_{\Omega} |\bigtriangledown u|^p dx+
  \int_{\Omega} g(u) u dx =\lambda \int_{\Omega} |u|^p dx \} .\]
The previous lemma says that for every $u\in W$, there exists  a
unique $t>0$ with $tu \in V$. We will show that if $u\in V$ , $u\geq
0$ and $F(u)=\min\limits_{u\in V} F(v)=\min\limits_{v\in W}
\max\limits_{s>0} F (sv)$ then $u$ is a solution for our problem.
$\square$

\begin{Lem} There exists $u\in V$ such that $F(u)=\min\limits_{v\in V}
F(v)$.
\end{Lem}
\paragraph{Proof.} Let $\{ u_n\}$ be a sequence in $V$ with
$F(u_n) \downarrow \inf\limits_{v\in V} F(v)$. Notice that we may
assume $u_n\geq 0$. We set $t_n= (\int_{\Omega} |\bigtriangledown
u_n|^p dx)^{1\over p}$ and $w_n= u_n/t_n$ for every $n\in N$. We
may assume that $\{w_n\}$ converges weakly in $V$ to some $w\in V$
and by Rellich theorem $\{w_n\}$ converges strongly to $w$ in $L^p
(\Omega)$. Moreover, by the Vitali convergence theorem
$\int_{\Omega} |w_n|^{1-\alpha} dx \to \int_{\Omega}
|w|^{1-\alpha} dx$.  We may assume $t_n \to t>0$, indeed, if $t_n
\to 0$, then we have
\[ \lambda \int_{\Omega} |w_n|^p dx= 1+ \int_{\Omega} {g(t_n w_n)
w_n^p \over t_n^{p} } dx\geq 1+ {m_1\over t_n^{\alpha+p-1}}
\int_{\Omega} |w_n|^{1-\alpha} dx \to  +\infty ,
\]
moreover if $t_n \to \infty$, then we must have
\[ \begin{array}{ll}
F(u_n) & \dis = \int_{\Omega} (\int_0^{u_n (x)} g(\tau) dt-
{1\over p} g (u_n) u_n) dx \\
& \dis \geq ( {m_1\over 1-\alpha} - {m_2\over p}) t_n^{1-\alpha}
\int_{\Omega} |w_n|^{1-\alpha} \to +\infty. \\
\end{array} \]
Thus a subsequence of $\{t_n\}$ converges to a positive number $
t$, then we have,
\[ 1+ \int_{\Omega} {g (tw)w \over t^{p-1}} dx= \lambda
\int_{\Omega} |w|^p dx. \]
 Now,  we will show $\int_{\Omega}
|\bigtriangledown w|^p dx=1$.
  Suppose not, then $\int_{\Omega} |\bigtriangledown w|^p dx<1$.
  By Lemma 2.2, there is $s\in (0,t)$ such that $sw\in V$. From
  (2.3) it follows that
    \[ \begin{array}{ll}
\dis  \inf_{v\in V} F(v) & \dis =\lim_{n\to \infty} F(u_n)=
\int_{\Omega} (\int_0^{t w(x)} g(\tau) dt - {1\over p} g(tw) tw)
dx \\
& \dis = \int_{\Omega} \int_0^{tw (x)} ((1-{1\over p}) g (\tau)-
1/p g' (\tau) \tau) dt dx \\
& \dis > \int_{\Omega} \int_0^{sw(x)} (1-{1\over p}) g(\tau) -
1/p g' (\tau) \tau dt dx= F(sw), \\
\end{array} \]
which is a contradiction. Therefore $\int_{\Omega}
|\bigtriangledown w|^p dx=1$, and hence $\{w_n\}$ converges
strongly to $w$ in $V$. This means that $tw\in V$ and
$F(tw)=\inf\limits_{v\in V} F(v)$.$\square$ \\

Now, we fix $u\in V$ with $F(u)=\min\limits_{v\in V} F(v)$. Since
$|u|\in V$ and $g$ is an odd function then $F(u)= F(|u|)$. Hence,
we may assume $u\geq 0$.

In this step, we show that $u>0$ in $\Omega$, which ensures
existence of the G\^ateaux derivative of $F$ at $u$ in the
direction of every $v\in C_0^{\infty} (\Omega) \cap U$.

\begin{Lem} If there is  $x_0 \in \Omega-\{0\}$ such that $u(x_0)=0$ then
$u_1 \equiv 0$, or $u_2 \equiv 0$, where
\[ u_1 (x) = \left\{ \begin{array}{ll}
u(x) & |x|\leq |x_0|, \\
0 & |x| \geq |x_0| \\
\end{array} \right.  \]
and
\[ u_2 (x) = \left\{ \begin{array}{ll}
0 & |x|\leq |x_0|, \\
u(x) & |x| \geq |x_0|. \\
\end{array} \right. \]
\end{Lem}
\paragraph{Proof.} Suppose that the conclusion does not hold,
i.e., there is $x_0 \in \Omega-\{0\}$ such that $u(x_0)=0$,
$u_1\not\equiv 0$ and $u_2\not\equiv 0$. From the definition of the
set $V$ we may assume $\int_{\Omega} |\bigtriangledown u_1|^p dx+
\int_{\Omega} g(u_1) u_1 dx \leq \lambda \int |u_1|^p dx $. By Lemma
2.2, there is $s\in (0,1]$ with $su_1 \in V$. Then (2.3) and
$u\not\equiv u_1$, implies  that $F(u)> F(s u_1)$, which is a
contradiction.$\square$

\begin{Lem} There is no $x_0 \in \Omega-\{0\}$ such that $u(x)=0$ for
every $x\in \Omega$ with $|x|\geq |x_0|$.
\end{Lem}

\paragraph{Proof.} Suppose that the conclusion does not hold.
Notice that  $u$ is not an eigenfunction of $-\bigtriangleup_p$
with homogeneous Dirichlet boundary condition for $\lambda_1$,
thus $\int_{\Omega} |\bigtriangledown u|^p >\lambda_1 \int
|u|^p$. Let $\epsilon$ be a  positive real number and
sufficiently small. For $s\in [1,1+\epsilon)$ we can define $u_s
\in W$ by $u_s (x) = u (x/s)$ for $x\in \Omega$. We set
\[ \varphi (t,s) = {t^p\over p} (s^{N-p} \int_{\Omega}
|\bigtriangledown u|^p dx- \lambda s^N  \int_{\Omega} |u|^p dx)
+s^N \int_{\Omega} \int_0^{tu (x)} g (\tau) \tau d\tau ,\]
 and
\[ \psi (t,s) =t^p (s^{N-p} \int_{\Omega} |\bigtriangledown u|^p
dx- \lambda s^N \int_{\Omega} |u|^p dx) + s^N \int_{\Omega} g(tu)
tu dx, \]
  for every $t,s \geq 0$. Notice that for $t>0$ and  $s\in [1,1+\epsilon)$
we will have  $\varphi (t,s) =F(t u_s)$ and
\[ \psi (t,s) =\int_{\Omega} |\bigtriangledown tu_s|^p
dx- \lambda \int_{\Omega} |tu_s|^p dx +  \int_{\Omega} g(tu_s)
tu_s dx. \]
 From $u\in V$ and (2.3), we obtain
 \[ \begin{array}{ll}
\dis  {\partial \psi \over \partial t} (1,1) & \dis
 = p (\int_{\Omega} |\bigtriangledown u|^p dx- \lambda
 \int_{\Omega} |u|^p)+\int_{\Omega} (g' (u) u^2 + u g(u)) dx \\
 & \dis = \int_{\Omega} g' (u) u^2 + (1-p) g (u) u dx<0 .\\
 \end{array} \]
Hence, the implicit function theorem implies that $\psi (t,s)=0$
defines a
 continuously differentiable function, $t= t(s)$ with $\psi
 (t(s), s)=0$ near $s=1$. On the other hand $\varphi (1,1)=\min\{
 \varphi (t(s),s)$ $1\leq s<1+\epsilon\}$, therefore
 \[ \begin{array}{l}
\dis 0\leq {\partial \varphi  \over \partial t} (1,1) {dt\over
ds}(1) + {\partial \varphi \over  \partial s} (1,1)={\partial
\varphi
\over \partial s} (1,1)\\

 \dis ={1\over p} (N-p) \int_{\Omega} |\bigtriangledown u|^p
 dx- {\lambda N\over p} \int_{\Omega} |u|^p dx+ N \int_{\Omega}
 \int_0^{u(x)} g(\tau) d\tau dx \\

\dis \leq \left( {N-p\over p} -{N m_2\over m_1 (1-\alpha)} \right)
\int_{\Omega} |\bigtriangledown u|^p dx+\lambda N \left(
{m_2\over m_1 (1-\alpha)}-{1\over p} \right) \int_{\Omega} |u|^p
dx \\

\dis < \left( {N-p \over p}-{N m_2\over m_1 (1-\alpha)} + {\lambda
N \over \lambda_1} ({m_2\over m_1 (1-\alpha)} -{1\over
p})\right) \int_{\Omega} |\bigtriangledown u|^p dx .\\
\end{array} \]
Thus, we must have
\[ \lambda_1 < \lambda - {\lambda p(1-\alpha) m_1\over p N m_2
- (N-p) (1-\alpha)m_1}, \]
 which contradicts (2.5). This completes the proof. $\square$
\begin{Lem} There is no $x_0 \in \Omega$  such that
 $u(x)=0$ for every $x\in \Omega$ with $|x| \leq |x_0|$.
 \end{Lem}

 \paragraph{Proof.} Let  $\Omega=\{x\in \mathbb{R}^N:
 |x| <R_1\}$. Suppose that the conclusion does not hold.
If $M$ is the maximum value of $u$ and $R$ is a point in $(0,
R_1)$ with $u(R) =M$.
 Then for $s\in [0,\epsilon)$,
  where $\epsilon$ is a sufficiently small positive real number,
 we can define $u_s\in W$ by
 \[ u_s (r) =\left\{ \begin{array}{ll}
 u(r+s) & 0\leq r\leq R-s, \\
 M & R-s\leq r \leq R, \\
 u(r) & R \leq r \leq R_1. \\
 \end{array} \right. \]
 Now, we define
 \[ \begin{array}{ll}
 \varphi (t,s)
& \dis = {t^p\over p} (\int_s^R |u'|^p (r-s)^{N-1} dr
+\int_R^{R_1} |u'|^p r^{N-1} dr \\
& \dis -\lambda \int_s^R |u|^p (r-s)^{N-1} dr- {\lambda M^p\over
N} (R^N- (R-s)^N) +\int_R^{R_1} |u|^p r^{N-1} dr)\\
& \dis +\int_s^R \int_0^{tu(r)} g(\tau) d\tau (r-s)^{N-1} dr + {R^N- (R-s)^N\over N} \int_0^{tM} g (\tau) d\tau \\
& \dis +\int_R^{R_1} \int_0^{tu(r)} g (\tau) d\tau r^{N-1} dr , \\
\end{array} \]
and
 \[ \begin{array}{ll}
 \psi (t,s)
 & \dis = t^p (\int_s^R |u'|^p (r-s)^{N-1} dr + \int_R^{R_1}
 |u'|^p r^{N-1} dr \\
& \dis -\lambda \int_s^R |u|^p (r-s)^{N-1} dr- {\lambda M^p\over
N} (R^N- (R-s)^N )\\
& \dis - \lambda \int_R^{R_1} |u|^p r^{N-1} dr) + \int_s^R g(tu)
tu (r-s)^{N-1} dr \\
 & \dis + {g (tM) tM\over N} (R^N- (R-s)^N) + \int_R^{R_1}
 g (tu) tu r^{N-1} dr . \\
\end{array} \]
Notice that $|S| \varphi (t,s)= F (tu_s)$ and
\[ |S| \psi (t,s) = \int_{\Omega} |\bigtriangledown (tu_s)|^p dx-
\lambda \int_{\Omega} |tu_s|^p dx+ \int_{\Omega} g (tu_s) tu_s dx
\]
for $t>0$ and $s\in [0,\epsilon)$, where $|S|$ is the measure of
the surface of the unit sphere $S$ in $\mathbb{R}^N$.\\
 From $u\in
V$, $\int_{\Omega} |\bigtriangledown u|^p dx> \lambda_1
\int_{\Omega} |u|^p dx$ and (2.4), we get
\[ \lambda \int_{\Omega} |u|^p= \int_{\Omega} |\bigtriangledown
u|^p dx+ \int_{\Omega} g(u) u dx >\lambda_1 \int_{\Omega} |u|^p
dx+ {m_1\over M^{p-1+\alpha}} \int |u|^p dx, \]
 which implies $\lambda -\lambda_1>{m_1\over M^{p-1+\alpha}}$.
 From ${\partial \psi \over \partial t} (1,0)<0$ and $\varphi
 (1,0)=\min \{ \varphi (t (s)),s): 0\leq s<\epsilon\}$, we obtain
 \[ \begin{array}{l}
\dis 0 \leq \underset{s\to 0^+}{\underline{\rm lim}} ({\partial
\varphi \over
\partial t} (t(s),s) {dt(s) \over ds} +{\partial \varphi\over
\partial s} (t(s),
s))=\underset{s\to 0^+}{\underline{\rm lim}}
 {\partial \varphi\over \partial s} (t(s),s) \\[0.4cm]
\dis ={1\over p} (- (N-1) \int_0^R |u'|^p r^{N-2} dr+ (N-1)
\lambda \int_0^R |u|^p r^{N-2} dr-\lambda M^p R^{N-1}) \\
\dis - (N-1) \int_0^R \int_0^{u(r)} g(\tau) d\tau r^{N-2} dr+
R^{N-1}
\int_0^M g(\tau) d\tau \\
\dis \leq {1\over p} ((N-1) \lambda \int_0^R |u|^p r^{N-2} dr
-\lambda M^p R^{N-1}) \\
\dis - (N-1) \int_0^R \int_0^{u(r)} g(\tau)d\tau r^{N-2} dr+
R^{N-1}
\int_0^M g(\tau) d\tau. \\
\end{array} \]
Since $\lambda -\lambda_1>{m_1\over M^{p-1+\alpha}} $ and
$H(\tau)={1\over \tau^p}\int_0^{\tau} g(p) dp$ is decreasing for
$\tau >0$, we have
\[ \begin{array}{ll}
- (N-1) &\dis \int_0^R |{u\over M}|^p r^{N-2} + R^{N-1} \\
        &\dis \leq {p\over \lambda M^P} (- (N-1) \int_0^R \int_0^{u(r)} g(\tau) d\tau r^{N-2}
        dr+ R^{N-1} \int_0^M g(\tau) d\tau) \\
        &\dis < {p (\lambda -\lambda_1) \over \lambda
        M^{1-\alpha}m_1} \int_0^M g(t) dt (- (N-1) \int_0^R
        |{u\over M}|^p r^{N-2} dr+ R^{N-1}). \\
\end{array} \]
Then, we obtain
\[ 1< {p (\lambda -\lambda_1) \over \lambda M^{1-\alpha} m_1} \int_0^M g(t) dt  \leq
      {p (\lambda -\lambda_1) \over \lambda M^{(1-\alpha)} m_1}{m_2\over 1-\alpha} M^{1-\alpha}=
      {p (\lambda -\lambda_1)m_2 \over \lambda (1-\alpha) m_1} , \]
or
\[ \lambda> {p \lambda_1 m_2 \over p m_2-(1-\alpha) m_1} . \]
 On the other hand
\[ \lambda -{\lambda (1-\alpha)m_1 \over p m_2} \leq \lambda -
{\lambda p (1-\alpha) m_1\over p N m_2- (1-\alpha) (N-p) m_1}\leq
\lambda_1.\] Therefore
\[ \lambda \leq  {p \lambda_1 m_2 \over p m_2-(1-\alpha) m_1} , \]
which is a contradiction.  This complete the proof. $\square$

\begin{Cor} For all  $x\in \Omega$, $u(x)\neq 0$.
\end{Cor}
 \paragraph{Proof.} It is a direct consequence of Lemmas 2.4, 2.5
 and 2.6. $\square$ \\

Now, we are ready to prove Theorem 2.1.

\paragraph{Proof of Theorem 2.1.}
By Corollary 2.6 we have $u>0$ on $\Omega$. Now, we will show
that $u$ is a weak solution of (1.1). In order to do this,we fix
$v \in c_0^{\infty} (\Omega) \cap U$ and define,
\[ \varphi (t,s)= {t^p\over p} (\int_{\Omega}  |\bigtriangledown (u+sv)|^p dx- \lambda \int_{\Omega}
 |u+sv|^p dx)+\]
\[ \int_{\Omega} \int_0^{t (u(x)+ sv(x))} g(t) dt dx ,\] \\
 and
\[ \psi (t,s)= t^p (\int_{\Omega} |\bigtriangledown (u+sv)|^p dx- \lambda \int_{\Omega}
|u+sv|^p dx )+\]
 \[\int_{\Omega} g(t(u+sv))t(u+sv)  dx ,\]\\
for  $t,s \in \mathbb{R}$. From $u\in V$ and (2.3), we have
${\partial \psi \over \partial t} (1,0)<0$. By  implicit function
theorem, $\psi (t,s)=0$ defines a continuously differentiable
function $t= t(s)$ with $\psi (t(s),s)=0$ near $s=0$.  Since for
some  $\epsilon >0$, $u\geq \epsilon$ on the support of $v$, the
function $F$ is G\^ateaux differentiable at $u$ in the direction
$v$.  This means that ${\partial \varphi \over
\partial s} (1,0)$ exists. Since $\varphi (1,0) =\min \{ \varphi
(t(s),s):s$ sufficiently close to $0\}$, we have
 \[ 0={\partial \varphi \over \partial t}(1,0) {dt\over ds} (0)+
 {\partial \varphi \over \partial s} (1,0) =\int_{\Omega}
 |\bigtriangledown u|^{p-2} \bigtriangledown u. \bigtriangledown v
 dx-\lambda \int_{\Omega} |u|^{p-2} u.v dx+ \]
 \[\int_{\Omega} g(u) v
 dx .\]\\
 Hence $u$ is a weak solution of problem (1.1). $\square$

 \section{Nonexistence result}
 \setcounter{equation}{5}
 Let  $\Omega \subset \mathbb{R}^N$ be a bounded
 domain with the  boundary $\partial \Omega$ which has the
 following  property: There exists a unit normal vector $v(x)=
 (v_1 (x), \ldots v_N (x))$ at every point $x\in \partial \Omega$
 and
 \begin{equation}
 \sum_{i=1}^N x_i v_i (x) \geq 0
 \end{equation}
for every $x= (x_1, \ldots, x_N) \in \partial \Omega$. Let us
consider the boundary value problem
\begin{equation}
\label{3.7}
\begin{array}{rl}
- \bigtriangleup_p u & = \lambda |u|^{p-2} u+ g(u) \quad
\text{in}\ \Omega,\\
 u & =0 \quad \text{on} \ \partial \Omega \\
\end{array}
\end{equation}
Here, we will show that this problem does not have a positive
solution in $W_0^{2, p_1} (\Omega) , (p_1>N)$. In oroder to see
this claim, let  $u\in W_0^{2,p_1} (\Omega), (p_1>N)$ be a
positive solution of this problem. By the Pohozaev identity
introduced in \cite{4}, we must have

\label{3.8}
\[ {N-p\over p} \int_{\Omega} |\bigtriangledown u|^p
dx- {\lambda N\over p} \int_{\Omega} |u|^p dx- N \int_{\Omega}
\int_0^{u(x)} g(t) dt= \]
            \[ - (1- {1\over p}) \int_{\partial
\Omega} |\bigtriangledown u|^p \sum_{i=1}^N x_i v_i dx.\]\\
On the other hand
\[ \int_{\Omega} |\bigtriangledown u|^p dx- \lambda \int_{\Omega}
|u|^p dx+ \int_{\Omega} g(u) u dx=0. \]
 From the above two identities, we see that the following identity holds for every $\beta \in
 \mathbb{R}$,
 \[ \begin{array}{ll}
 \dis ({N-p \over p}+\beta)
& \dis  \int_{\Omega} |\bigtriangledown u|^p d x- \lambda ({N\over
p} +\beta)\int_{\Omega} |u|^p dx - N \int_{\Omega} \int_0^{u(x)} g(t) dt + \\
 & \dis \beta \int_{\Omega} g(u) u  dx= - (1-{1\over p})
 \int_{\partial \Omega} |\bigtriangledown u|^p \sum_{i=1}^N x_i v_i ds. \\
\end{array} \]
Hence
\begin{equation}
\label{3.9}
\begin{array}{ll}
\dis   ({N-p\over p}+\beta)
  & \dis \int_{\Omega} |\bigtriangledown u|^p dx -
  \lambda ({N\over p}+\beta) \int_{\Omega} |u|^p dx- ({Nm_2\over
  1-\alpha}+\beta m_1)\\
& \dis  \int_{\Omega} u^{1-\alpha} dx \leq - (1-{1\over  p})
 \int_{\partial \Omega} |\bigtriangledown u|^p
\sum_{i=1}^N   x_i v_i(x) ds .\\
\end{array}
  \end{equation}
  Now, it follows from (\ref{3.9}) that the following inequalities
\begin{equation}
\label{3.10} {N-p\over p} +\beta \geq 0,
\end{equation}
\begin{equation}
\label{3.11}
 - \lambda ({N\over p} +\beta) \geq 0,
\end{equation}
\begin{equation}
\label{3.12}
  -( {N m_2\over 1-\alpha} +\beta m_1) \geq 0,
\end{equation}
cannot hold simultaneously with at least one strict inequality sign.
Thus, we have the following nonexistence result.

\begin{Th}Let $N\geq 2$, $p>1$, $0<\alpha<1$, $m_1>0$, $m_2>0$ be real
numbers such that
 (\ref{3.10}),  (\ref{3.11}) and (\ref{3.12})
hold with at least one strict inequality sign. Then the
boundary-value problem (\ref{3.7}), has no positive solution  in
$W^{2, p_1} (\Omega)\cap W_0^{1, p_1}(\Omega)$ for $p_1>N$.
\end{Th}

%

\end{document}